\documentclass[11pt]{amsart}
\usepackage[T1]{fontenc}
\usepackage[utf8]{inputenc}
\usepackage[margin=0.95in]{geometry}
\usepackage{amsmath,amssymb,amsthm,mathtools}
\usepackage{newtxtext,newtxmath}
\usepackage{microtype,enumitem}
\usepackage[hidelinks]{hyperref}
\usepackage{comment}
\usepackage{xcolor}

\AddToHook{env/thebibliography/after}{\enlargethispage{2\baselineskip}}
\setlist[enumerate]{leftmargin=*,itemsep=0.4em,topsep=0.4em}
\newtheorem{theorem}{Theorem}[section]
\newtheorem{lemma}[theorem]{Lemma}
\newtheorem{proposition}[theorem]{Proposition}
\newtheorem{corollary}[theorem]{Corollary}
\newtheorem{maintheorem}{Theorem}

\theoremstyle{definition}

\theoremstyle{remark}
\newtheorem{remark}[theorem]{Remark}
\newcommand{\D}{\mathbb D}
\newcommand{\C}{\mathbb C}
\newcommand{\R}{\mathbb R}
\newcommand{\Z}{\mathbb Z}

\DeclareMathOperator{\Area}{Area}

\DeclareMathOperator{\wind}{wind}
\DeclareMathOperator{\rank}{rank}
\DeclareMathOperator{\tr}{tr}

\DeclareMathOperator{\Ric}{Ric}

\title[Holomorphicity with Lagrangian free boundary]{Holomorphicity of minimal surfaces\\with Lagrangian free boundary}
\author{Jingyi Chen}
\date{October 16, 2026}
\subjclass[2020]{Primary 58E12; Secondary 53D12, 53C21}
\thanks{The research is partially supported by NSERC Discovery Grant GR010074.}
\address{Department of Mathematics, The University of British Columbia, Vancouver, B.C. V6T 1Z2, Canada}
\email{jychen@math.ubc.ca}
\begin{document}
\begin{abstract}
We study branched minimal immersions of compact bordered Riemann
surfaces into almost K\"ahler four-manifolds with free boundary on
oriented Lagrangian surfaces. We establish a relative complex-point
formula relating the Maslov index to complex and anticomplex
tangency counts, and derive a sharp criterion for holomorphicity
and area minimization in the relative homotopy class. If the
ambient manifold is K\"ahler--Einstein with positive Einstein
constant and the supporting Lagrangian surface is minimal, then
the immersion is holomorphic, antiholomorphic, or Lagrangian
whenever $\cos\alpha$ does not change sign, where $\alpha$ is the
K\"ahler angle. In the hyperk\"ahler case, under the same boundary
hypotheses, the absence of Lagrangian points implies holomorphicity
with respect to one fixed parallel complex structure. We also
extend Lichnerowicz's partial-energy homotopy identity to compact
almost K\"ahler domains with boundary mapped into a Lagrangian
submanifold. Consequently, holomorphic and antiholomorphic maps
minimize energy in their relative homotopy classes. Examples show
that total reality cannot replace the Lagrangian boundary condition.
\end{abstract}
\maketitle

\section{Introduction}\label{sec:introduction}
Complex and anticomplex tangencies of minimal surfaces in K\"ahler
four-manifolds are related to characteristic classes. For closed
surfaces, Webster's formulas express their algebraic counts in terms
of the first Chern class and related topological data
\cite{Webster,WebsterChern}. Chen--Tian extended this study to almost
K\"ahler four-manifolds~\cite{ChenTian}.

We study smooth nonconstant branched minimal immersions
$$
u:(\Sigma,\partial\Sigma)\longrightarrow(M^4,L^2)
$$
from a compact connected bordered Riemann surface, satisfying the
free-boundary condition, where $(M,g,J,\omega)$ is almost K\"ahler
and $L$ is an oriented embedded Lagrangian surface.
Let $P(u)$ and $Q(u)$ be the complex and anticomplex
tangency counts, including tangencies at branch points.  Interior
multiplicities are counted in full and boundary multiplicities with
one half of the reflected order.  We write $\mu_L(u)$ for the Maslov
index of the bundle pair $(u^*TM,u|_{\partial\Sigma}^*TL)$.
We also put $\omega(u)=\int_\Sigma u^*\omega$ and write
$\Area_g(u)$ for parametrized area.

Our first result is a relative complex-point identity and its
holomorphicity consequence.

\begin{maintheorem}\label{thm:intro-A}
Let $(M^4,g,J,\omega)$ be almost K\"ahler and let
$L^2\subset M$ be a smooth oriented embedded Lagrangian surface.
Let $u:(\Sigma,\partial\Sigma)\to(M,L)$ be a smooth nonconstant
branched minimal immersion satisfying the free-boundary condition.
Then the following hold.
\begin{enumerate}[label=(\roman*)]
\item If $u$ is neither $J$-holomorphic nor $J$-antiholomorphic, then
$$
\mu_L(u)=2\bigl(Q(u)-P(u)\bigr).
$$
\item If $u$ is not $J$-antiholomorphic and
$$
\mu_L(u)>2Q(u),
$$
then $u$ is $J$-holomorphic and minimizes area in its relative
homotopy class.  In particular, this holds if $u$ has no
anticomplex points and $\mu_L(u)>0$.
\end{enumerate}
\end{maintheorem}

The ramification and Euler-characteristic contributions cancel
in the derivation of part~(i).
The Lagrangian hypothesis cannot in general be replaced by total
reality, and the strict inequality in part~(ii) is sharp, even among
immersed area-minimizing disks; see Section~\ref{sec:examples-sharpness}.

We also extend Lichnerowicz's partial-energy homotopy identity
from closed almost K\"ahler domains~\cite{Lichnerowicz} to compact
domains with boundary mapped into a Lagrangian submanifold.
The boundary term in the homotopy variation vanishes under this
constraint. Consequently, holomorphic and antiholomorphic maps
minimize energy in their relative homotopy classes. For bordered
Riemann surfaces, Wirtinger's inequality also gives area minimization
in the relative homotopy class and free-boundary
stationarity~\cite{HarveyLawson}.

Stability provides another approach to holomorphicity.
Siu--Yau~\cite{SiuYau} proved holomorphicity for stable harmonic two-spheres under positive
holomorphic bisectional curvature, a key ingredient in their proof
of the Frankel conjecture. Related stability-based results were
obtained by Micallef~\cite{Micallef}, Micallef--Wolfson~\cite{MW},
Arezzo~\cite{Arezzo}, and Arezzo--Micallef~\cite{ArezzoMicallef}. A complementary
approach based on moving frames and the K\"ahler angle was developed
by Chern--Wolfson~\cite{ChernWolfson} and Wolfson~\cite{Wolfson}.
In the free-boundary setting, Chen--Fraser~\cite{CF} studied admissible
holomorphic variations of minimal disks with Lagrangian
boundary. They also proved holomorphicity and rigidity
results for stable free-boundary minimal disks on weakly two-convex
hypersurfaces under a nonnegative isotropic curvature assumption,
including K\"ahler and hyperk\"ahler cases~\cite{CFCrelle}.

We use $\pm J$-holomorphic to mean either $J$-holomorphic or
$J$-antiholomorphic. On the regular set of $u$, write
$u^*\omega=\cos\alpha\,dA_u$, where $dA_u$ is the area form of
$u^*g$ and $\alpha$ is the K\"ahler angle.
The function $\cos\alpha$ extends continuously across the branch
points of $u$.

For a K\"ahler--Einstein four-manifold, equation~\eqref{eq:maslov-integral}
relates the Maslov index to the symplectic area and the mean curvature
of the supporting Lagrangian surface $L$~\cite{CG,Pacini}.

\begin{maintheorem}\label{thm:intro-B}
Let $(M^4,g,J,\omega)$ be K\"ahler--Einstein with
$\Ric_M=\kappa g$, $\kappa>0$, and let $L^2\subset M$ be an oriented
minimal Lagrangian surface.  Let
$u:(\Sigma,\partial\Sigma)\to(M,L)$ be a smooth nonconstant branched
minimal immersion satisfying the free-boundary condition.
Then the following hold.
\begin{enumerate}[label=(\roman*)]
\item For each $\varepsilon\in\{+1,-1\}$, if
$\varepsilon\cos\alpha\ge0$ on $\Sigma$, then $u$ is either
$\varepsilon J$-holomorphic and minimizes area in its relative
homotopy class, or Lagrangian. In particular, if $u$ has no
Lagrangian points, then it is $J$-holomorphic or
$J$-antiholomorphic and minimizes area in its relative homotopy class.
\item For every $A\in H_2(M,L;\Z)$,
$$
\mu_L(A)=\frac{\kappa}{\pi}\,\omega(A).
$$
\end{enumerate}
\end{maintheorem}

For a map $u$ in Theorem~B that is not $\pm J$-holomorphic, Theorem~A and
\eqref{eq:maslov-integral} yield
$$
\omega(u)=\frac{2\pi}{\kappa}\bigl(Q(u)-P(u)\bigr).
$$
Thus zero symplectic area is equivalent to $P(u)=Q(u)$ under these
hypotheses; in particular, a totally real branched minimal surface
has zero symplectic area. The proof of rigidity when
$\cos\alpha$ has one sign uses the complex-point formula rather than
the K\"ahler-angle differential equation.  This is the relative counterpart of Webster's approach
for closed minimal surfaces~\cite{Webster,WebsterChern}.

When the Einstein constant vanishes and $L$ is minimal Lagrangian,
equation~\eqref{eq:maslov-integral} gives $\mu_L(u)=0$.
Theorem~\ref{thm:intro-A}(i) then implies $P(u)=Q(u)$ unless $u$ is
$\pm J$-holomorphic. In the hyperk\"ahler case, the K\"ahler-angle
equation yields a further rigidity statement:

\begin{maintheorem}\label{thm:intro-C}
Let $(M^4,g,J,\omega)$ be hyperk\"ahler, and let
$L^2\subset M$ be an oriented minimal Lagrangian surface. Let
$u:(\Sigma,\partial\Sigma)\to(M,L)$ be a smooth nonconstant
branched minimal immersion satisfying the free-boundary condition.
If $u$ has no Lagrangian points, then it is holomorphic with
respect to one fixed parallel complex structure in the
hyperk\"ahler $S^2$-family.
\end{maintheorem}

The proof of Theorem~C combines the K\"ahler-angle equation and its
Neumann boundary condition with the hyperk\"ahler Gauss-map argument;
compare Wolfson's closed-surface rigidity results~\cite{Wolfson}.

Section~\ref{sec:AK-framework} develops the bundle-pair conventions,
complex tangency orders, and boundary factorization for general
bordered surfaces.  Section~\ref{Lichnerowicz} proves the Lagrangian-boundary
extension of Lichnerowicz's homotopy identity and its minimizing
consequences.  Section~\ref{sec:rigidity} proves the rigidity results and
their K\"ahler--Einstein refinements.  Section~\ref{sec:examples-sharpness}
provides sharpness examples and counterexamples to possible
strengthenings.

\medskip
\noindent\textit{Acknowledgement.}
This work was carried out in part during the author's visit to BICMR
in 2025--26.  He thanks Professor Gang Tian and the institute for
their hospitality and support, and Professor Julien Duval for
helpful communications concerning an earlier version of this work.
The author acknowledges the use of ChatGPT (OpenAI) for assistance
in checking some examples in
Section~\ref{sec:examples-sharpness}, as well as for editorial
suggestions during the preparation of the manuscript.

\section{Maslov index, complex tangencies, and boundary problems}
\label{sec:AK-framework}
Throughout Sections~2--4, $(\Sigma,j)$ is a compact connected bordered
Riemann surface with nonempty smooth boundary and the orientation determined by $j$.
If $g_\Sigma$ is the genus and $b\ge1$ the number of boundary
components, then the Euler characteristic of $\Sigma$ is $\chi(\Sigma)=2-2g_\Sigma-b$.
For the minimal surface results in Sections~2 and~4, the ambient
manifold $(M,g,J,\omega)$ is almost K\"ahler, unless stronger
hypotheses are stated, that is, $J$ is an almost complex structure
satisfying $g(JX,JY)=g(X,Y)$, $\omega(X,Y)=g(JX,Y)$, and $d\omega=0$.

\subsection{Bundle pair and Maslov index}\label{subsec:AK-Maslov}
We recall the standard bundle-pair description of the Maslov index for
an oriented totally real boundary condition, fixing the conventions
used below; see~\cite[Appendix~C]{MS} and~\cite{Pacini}.

Let $L^2\subset M$ be a smooth oriented embedded totally real surface,
so that $T_xL\cap JT_xL=\{0\}$ for every $x\in L$, and let
$u:(\Sigma,\partial\Sigma)\to(M,L)$ be smooth.
The complexification of the pullback tangent bundle is
$\mathcal E=u^*TM\otimes_\R\C$.  Complex conjugation is defined by
$\overline{X+iY}=X-iY$ for real sections $X,Y$ of $u^*TM$.
The $i$ and $-i$ eigenspaces of $J$ are
$$
E=u^*T^{1,0}M=\{V\in\mathcal E:JV=iV\},\quad
\overline E=u^*T^{0,1}M=\{V\in\mathcal E:JV=-iV\}.
$$
Thus $\mathcal E=E\oplus\overline E$, with
$\rank_\C\mathcal E=4$ and $\rank_\C E=\rank_\C\overline E=2$.
The tangent bundle of $L$ determines the bundle pair
\begin{equation}\label{eq:EF}
(E,F),\quad
F_z=\{X-iJX:X\in T_{u(z)}L\}\subset E_z
\quad(z\in\partial\Sigma).
\end{equation}
The map $X\mapsto X-iJX$ identifies multiplication by $J$ with
multiplication by $i$, since $JX-iJ(JX)=i(X-iJX)$.  Thus total reality
of $L$ gives
$$
\rank_\R F=2,\quad F\cap iF=\{0\},\quad E|_{\partial\Sigma}=F\oplus_\R iF.
$$
The conjugate boundary bundle is
$\overline F_z=\{X+iJX:X\in T_{u(z)}L\}\subset\overline E_z$.

Write $g$ also for the pullback metric on $u^*TM$, and
$(\cdot,\cdot)$ for its complex-bilinear extension to $\mathcal E$.
The associated Hermitian inner product
$\langle V,W\rangle=(V,\overline W)$ is linear in the first argument.
The subbundles $E$ and $\overline E$
are isotropic for $(\cdot,\cdot)$, and the pairing between
$E$ and $\overline E$ is nondegenerate.

For real $X,Y\in TL$, direct expansion gives
\begin{equation}\label{eq:metric-types}
\langle X-iJX,Y-iJY\rangle
=2g(X,Y)-2i\omega(X,Y).
\end{equation}
When $L$ is Lagrangian, the restriction of the
Hermitian inner product to $F$ is real and positive
definite.

Let $\mu(E,F)$ be the Maslov index of the bundle pair. We recall
its homological invariance and parity for an oriented totally real
boundary condition; see~\cite[Appendix~C]{MS}.

\begin{proposition}\label{prop:class-positivity}
Let $L^2\subset M^4$ be a smooth oriented totally real surface.  Then, for
any smooth map $u:(\Sigma,\partial\Sigma)\to(M,L)$, the Maslov index
of its bundle pair satisfies
$$
\mu_L(u):=\mu(E,F)
=\bigl\langle\mathfrak m_L,u_*[\Sigma,\partial\Sigma]\bigr\rangle
\in2\Z,
$$
where $\mathfrak m_L\in H^2(M,L;\Z)$ is the relative Maslov class, and 
$\mu_L(u)$ depends only on the relative homology class.
\end{proposition}

\begin{proof}
Put $D=\det_\C T^{1,0}M$. The orientation of $L$ determines a
unit section $\tau_L$ of $D|_L$: for any positively oriented frame
$X_1,X_2$ of $TL$, normalize
$(X_1-iJX_1)\wedge(X_2-iJX_2)$ to have unit length.
Total reality makes this wedge nonzero, and a change of oriented
frame multiplies it by a positive real number. Thus $\tau_L$ is
well defined, and we set
$\mathfrak m_L=2c_1(D,\tau_L)\in H^2(M,L;\Z)$.

Since $\Sigma$ has boundary, $E$ admits a global complex
trivialization.  On each oriented boundary component $C\subset\partial\Sigma$,
let $U_C:C\to GL(2,\C)$ be the matrix of a positively oriented
real frame of $F$ in this trivialization.  Then
\begin{equation}\label{eq:Maslov-winding}
\mu(E,F)=\sum_{C\subset\partial\Sigma}
\wind_C\left(\frac{(\det_\C U_C)^2}{|\det_\C U_C|^2}\right)
=2\sum_{C\subset\partial\Sigma}\wind_C(\det_\C U_C).
\end{equation}
The orientations are the boundary orientations induced by $j$.
The sum of winding numbers is the relative first Chern number of
$(u^*D,u|_{\partial\Sigma}^*\tau_L)$. Naturality therefore identifies
\eqref{eq:Maslov-winding} with
$\langle\mathfrak m_L,u_*[\Sigma,\partial\Sigma]\rangle$,
proving both homological invariance and evenness.
\end{proof}

For $A\in H_2(M,L;\Z)$ we also write
$\mu_L(A)=\langle\mathfrak m_L,A\rangle$.  When $L$ is Lagrangian,
we use
$$
\omega(u)=\int_\Sigma u^*\omega,\quad
\omega(A)=\langle[\omega],A\rangle,
$$
where the latter is a relative pairing because $\omega|_L=0$.

\subsection{Differential components and complex tangencies}
\label{subsec:AK-complex-points}
Let $u:\Sigma\to M$ be a smooth nonconstant weakly conformal
harmonic map, which we call a branched minimal immersion. Write
$K_\Sigma=(T^{1,0}\Sigma)^*$ for the canonical line bundle. 

In a local holomorphic coordinate $z=x+iy$, subscripts $x,y,z,\bar z$ on $u$
or $\nabla$ denote differentiation in this coordinate. 
Write $u_x=du(\partial_x)$, $u_y=du(\partial_y)$ and
$$
u_z=\tfrac12(u_x-iu_y),\quad
u_{\bar z}=\tfrac12(u_x+iu_y)=\overline{u_z}.
$$
In this coordinate define the real type components
\begin{equation}\label{eq:partial-J-definitions}
\partial_Ju=\tfrac12(u_x-Ju_y),
\quad
\bar\partial_Ju=\tfrac12(u_x+Ju_y),
\end{equation}
and their complex representatives
\begin{equation}\label{eq:Phi-Psi}
\begin{aligned}
\Phi&=\partial_Ju-iJ\partial_Ju=2(u_z)^{1,0}\in E,\\
\Psi&=\bar\partial_Ju+iJ\bar\partial_Ju
=2(u_z)^{0,1}\in\overline E,
\end{aligned}
\end{equation}
so that $u_z=\tfrac12(\Phi+\Psi)$.
On each holomorphic coordinate neighborhood $U\subset\Sigma$,
$\Phi$ and $\Psi$ are local sections of $E$ and
$\overline E$, respectively. Writing $\Phi(z)$ and $\Psi(z)$
to indicate the chosen coordinate, a holomorphic change
$z=z(w)$ gives
$$
\Phi(w)=\frac{dz}{dw}\Phi(z),\quad
\Psi(w)=\frac{dz}{dw}\Psi(z).
$$
Thus the globally defined bundle-valued $(1,0)$-forms are
\begin{equation}\label{eq:Phi-Psi-global}
\boldsymbol\Phi=\Phi\,dz\in\Gamma(K_\Sigma\otimes E),\quad
\boldsymbol\Psi=\Psi\,dz\in\Gamma(K_\Sigma\otimes\overline E).
\end{equation}

At a regular point $z\in\Sigma$ of $u$ is 
a \emph{complex point} if $du\circ j=J\circ du$ on $T_z\Sigma$, and an
\emph{anticomplex point} if $du\circ j=-J\circ du$ on $T_z\Sigma$.
By \eqref{eq:Phi-Psi}, these conditions are equivalent
to $\Psi(z)=0$ and $\Phi(z)=0$, respectively.
A regular point $z$ is \emph{totally real} if
$du(T_z\Sigma)\cap Jdu(T_z\Sigma)=\{0\}$.
At branch points, these three notions are defined using the extended
tangent plane constructed in Section~\ref{subsec:tangency-counts}.
We call $u$ \emph{totally real} if it is totally real at every point
of $\Sigma$, including the boundary.

On the regular set, in a local chart $z=x+iy$ write
\begin{equation}\label{eq:Kahler-angle}
 u^*g=\lambda^2(dx^2+dy^2),
 \quad
 dA_u=\lambda^2\,dx\wedge dy,
 \quad
 u^*\omega=\cos\alpha\,dA_u,
\end{equation}
where $\lambda>0$ and $\alpha\in[0,\pi]$ is the K\"ahler angle.
A regular point is complex, Lagrangian, or anticomplex precisely
when $\cos\alpha$ equals $1$, $0$, or $-1$, respectively;
it is totally real precisely when $|\cos\alpha|<1$.
Since $g(u_x,Ju_y)=-\lambda^2\cos\alpha$, we obtain
\begin{equation}\label{eq:partial-angle}
\begin{aligned}
|\partial_Ju|_g^2
&=\tfrac14\bigl(|u_x|_g^2+|u_y|_g^2-2g(u_x,Ju_y)\bigr)
=\tfrac12\lambda^2(1+\cos\alpha),\\
|\bar\partial_Ju|_g^2
&=\tfrac14\bigl(|u_x|_g^2+|u_y|_g^2+2g(u_x,Ju_y)\bigr)
=\tfrac12\lambda^2(1-\cos\alpha).
\end{aligned}
\end{equation}
The corresponding Hermitian norms are
$|\Phi|^2=2|\partial_Ju|_g^2$ and
$|\Psi|^2=2|\bar\partial_Ju|_g^2$.

Let $\nabla$ denote the pullback of the Levi--Civita connection of $g$,
extended complex-linearly to $\mathcal E$.  Thus
$\nabla_x=\nabla_{\partial_x}$, $\nabla_y=\nabla_{\partial_y}$ and
$\nabla_{\bar z}=\tfrac12(\nabla_x+i\nabla_y)$ on pullback sections.
Define
\begin{equation}\label{eq:dbar-operators}
\begin{aligned}
\bar\partial^\nabla s&=\nabla_{\bar z}s\,d\bar z,
&&s\in\Gamma(\mathcal E),\\
\bar\partial_E s&=(\nabla_{\bar z}s)^{1,0}d\bar z,
&&s\in\Gamma(E),\\
\bar\partial_{\overline E}t&=(\nabla_{\bar z}t)^{0,1}d\bar z,
&&t\in\Gamma(\overline E).
\end{aligned}
\end{equation}
These are complex-linear $(0,1)$-connections on
$\mathcal E$, $E$, and $\overline E$, respectively; for example,
$\bar\partial_E(fs)=(\bar\partial f)s+f\bar\partial_Es$. 
Let $\pi_+=\tfrac12(I-iJ)$ and $\pi_-=\tfrac12(I+iJ)$ be the
projections of $\mathcal E$ onto $E$ and $\overline E$, respectively. 
In the almost K\"ahler case, $\nabla$ need not preserve $E$ and
$\overline E$, so the projections in \eqref{eq:dbar-operators}
are essential. Since the base is a Riemann surface, each operator
defines a holomorphic bundle structure by the Koszul--Malgrange
theorem~\cite{KM}. A local holomorphic frame $e_1,e_2$ of $E$ thus
satisfies $\bar\partial_Ee_1=\bar\partial_Ee_2=0$.

Since the Levi--Civita connection of $g$ is torsion-free and $u$
is harmonic, we have $\nabla_xu_y=\nabla_yu_x$ and
$\nabla_xu_x+\nabla_yu_y=0$. Therefore
\begin{equation}\label{eq:MM-uz}
\nabla_{\bar z}u_z
=\tfrac14(\nabla_xu_x+\nabla_yu_y)
 +\tfrac{i}{4}(\nabla_yu_x-\nabla_xu_y)=0.
\end{equation}
Thus $u_z$ is a local holomorphic section of $\mathcal E$
with respect to $\bar\partial^\nabla$, as observed by
Micallef--Moore~\cite{MM}. Under a holomorphic coordinate
change $z=z(w)$,
$$
u_w=\frac{dz}{dw}u_z,\quad u_w\,dw=u_z\,dz.
$$
Consequently, the local forms $u_z\,dz$ define a global
holomorphic section of $K_\Sigma\otimes\mathcal E$, with
its tensor-product holomorphic structure. In the notation of
\eqref{eq:Phi-Psi-global}, this section is
$\tfrac12(\boldsymbol\Phi+\boldsymbol\Psi)$.

Choose a local $\bar\partial^\nabla$-holomorphic frame
$e_1,\ldots,e_4$ of $\mathcal E$. In this frame, 
$$
u_z=\sum_{j=1}^4 f^je_j,\quad
0=\nabla_{\bar z}u_z=\sum_{j=1}^4(\partial_{\bar z}f^j)e_j,
$$
so the coefficient functions $f^j$ are holomorphic. At an interior
branch point $z_0$, let $r\ge1$ be the least vanishing order among
the nontrivial $f^j$. This is the branching order of $u$, and
\begin{equation}\label{eq:LC-ramification}
u_z=(z-z_0)^rV,\quad
\nabla_{\bar z}V=0,\quad V(z_0)\ne0.
\end{equation}

\subsection{Boundary conditions and local factorization}
\label{subsec:boundary-factorization}

\begin{lemma}\label{lem:boundary-lagrangian}
Let $L^2\subset M^4$ be Lagrangian and let
$u:(\Sigma,\partial\Sigma)\to(M,L)$ be a smooth
weakly conformal map satisfying the free-boundary
condition. In a boundary-adapted holomorphic
coordinate $z=s+it$, with $\partial\Sigma$ given
by $t=0$ and the interior by $t>0$, one has
\begin{equation}\label{eq:boundary-components}
\Phi\in F,\quad \Psi\in\overline F
\quad\text{along }t=0.
\end{equation}
\end{lemma}

\begin{proof}
Along $t=0$, the free-boundary condition gives $u_s\in TL$ and
$u_t\in(TL)^\perp=JTL$. Consequently,
$$
\partial_Ju=\tfrac12(u_s-Ju_t)\in TL,
\quad
\bar\partial_Ju=\tfrac12(u_s+Ju_t)\in TL.
$$
The definitions of $\Phi$ and $\Psi$ give
\eqref{eq:boundary-components}.
\end{proof}

\begin{lemma}\label{lem:relative-factorization}
Let $L^2\subset M^4$ be Lagrangian, and let
$u:(\Sigma,\partial\Sigma)\to(M,L)$ be a smooth
nonconstant branched minimal immersion satisfying
the free-boundary condition. 
If $\Phi\not\equiv0$, then its zero set in $\Sigma$ is
finite. At an interior zero $z_0$,
\begin{equation}\label{eq:Phi-interior-factor}
\Phi=(z-z_0)^m e^\sigma s,
\quad m\ge1,\quad
\bar\partial_Es=0,\quad s(z_0)\ne0,
\end{equation}
where $\sigma\in W^{1,p}_{\mathrm{loc}}$ for every
$p>2$. At a boundary zero, in a boundary-adapted
holomorphic coordinate $\zeta$ centered at the zero 
with $\partial\Sigma$ represented by the real axis, 
\begin{equation}\label{eq:Phi-boundary-factor}
\Phi=\zeta^m V,
\quad m\ge1,\quad V(0)\ne0,
\end{equation}
where $V\in W^{1,p}$ for every $p>2$, and $V$ has boundary values in
$F$ along the real axis. 
Moreover, the complex line subbundle $K_\Phi=\mathbb C\Phi\subset E$ on
$\{\Phi\ne0\}$ extends continuously across the zeros of $\Phi$, and $F_\Phi=K_\Phi|_{\partial\Sigma}\cap F$ is a real line
subbundle along $\partial\Sigma$. 
The analogous conclusions hold for $\Psi$, with
$K_\Psi\subset\overline E$ and
$F_\Psi=K_\Psi|_{\partial\Sigma}\cap\overline F$.
\end{lemma}

\begin{proof}
{\bf Step 1.} \emph{Component equations.}
Differentiating $J^2=-I$ gives
$$
J(\nabla_AJ)=-(\nabla_AJ)J.
$$
Then, for every complexified tangent vector $A$,
$$
(\nabla_AJ)(T^{1,0}M)\subset T^{0,1}M,
\quad
(\nabla_AJ)(T^{0,1}M)\subset T^{1,0}M.
$$
The identity $d\omega=0$ now gives
\begin{equation}\label{eq:mixed-type}
(\nabla_{\overline X}J)Y=0,
\quad
(\nabla_XJ)\overline Y=0
\quad\mbox{for $X,Y\in T^{1,0}M$}.
\end{equation}
Indeed, for $X,Y,Z\in T^{1,0}M$, the cyclic identity for $d\omega$
gives
$$
((\nabla_{\overline X}J)Y,Z)
+((\nabla_YJ)Z,\overline X)
+((\nabla_ZJ)\overline X,Y)=0.
$$
The complex-bilinear metric vanishes on
$T^{1,0}M\times T^{1,0}M$ and on $T^{0,1}M\times T^{0,1}M$.
Since $(\nabla_YJ)Z\in T^{0,1}M$ and
$(\nabla_ZJ)\overline X\in T^{1,0}M$, the second and third terms
pair vectors of the same type and hence vanish. Since
$(\nabla_{\overline X}J)Y\in T^{0,1}M$ and the pairing between
$T^{0,1}M$ and $T^{1,0}M$ is nondegenerate, the first identity follows.
The second follows in the same way from the cyclic identity with
$X,\overline Y,\overline Z$.

From $\Phi=2\pi_+u_z$, $\Psi=2\pi_-u_z$, and
\eqref{eq:MM-uz},
\begin{equation}\label{eq:component-equations}
\begin{aligned}
(\nabla_{\bar z}\Phi)^{1,0}
 &=-\tfrac{i}{4}(\nabla_{\overline\Phi}J)\Psi,\\
(\nabla_{\bar z}\Psi)^{0,1}
 &=\tfrac{i}{4}(\nabla_{\overline\Psi}J)\Phi.
\end{aligned}
\end{equation}
The complex-bilinear metric vanishes on \(E\times E\) and \(\overline E\times\overline E\), so \((\Phi,\Phi)=(\Psi,\Psi)=0\). Expanding \(u_z=\tfrac12(\Phi+\Psi)\) and using conformality, we obtain
\[
(\Phi,\Psi)=2(u_z,u_z)=0.
\]
Set $W=(\nabla_{\overline\Phi}J)\Psi$. Since $J\Psi=-i\Psi$
and $J(\nabla_{\overline\Phi}J)=-(\nabla_{\overline\Phi}J)J$,
we have $JW=iW$, so $W\in E$. Since $\nabla J$ is skew-adjoint
with respect to the complex-bilinear pairing,
$$
(W,\Psi)=((\nabla_{\overline\Phi}J)\Psi,\Psi)=0.
$$
If $\Psi\ne0$, the kernel in $E$ of the functional
$Z\longmapsto(Z,\Psi)$ is a one-dimensional complex subspace.
It contains $\Phi$, and hence equals $\C\Phi$ wherever
$\Phi\ne0$.  If $\Psi=0$, then $W=0$.  Thus the
right-hand side of the first equation in
\eqref{eq:component-equations} is a scalar multiple of $\Phi$ wherever
$\Phi\ne0$, and it vanishes where $\Phi=0$.  

Writing $W=a_0\Phi$ wherever $\Phi\ne0$, we have
$$
|a_0|\,|\Phi|=|W|
\le|\nabla J|\,|\Phi|\,|\Psi|,
$$
so $|a_0|\le c_0|\Psi|$. Set $a_0=0$ at
the zeros of $\Phi$. The same argument, with
$\Phi$ and $\Psi$ interchanged, applies to the
second equation. Consequently, there exist
$a,b\in L^\infty_{\mathrm{loc}}$ such that
\begin{equation}\label{eq:scalar-component-equations}
\bar\partial_E\Phi=a\Phi\,d\bar z,
\quad
\bar\partial_{\overline E}\Psi=b\Psi\,d\bar z,
\end{equation}
with $|a|\le c_0|\Psi|$ and $|b|\le c_0|\Phi|$; compare
\cite[eqs.~(46), (49), pp.~901--903]{ChenTian}.

{\bf Step 2.} \emph{Interior and boundary factorization.}

\emph{Interior.}
In a local $\bar\partial_E$-holomorphic frame $e_1,e_2$, write
$\Phi=f^1e_1+f^2e_2$. Equation~\eqref{eq:scalar-component-equations}
then gives $\partial_{\bar z}f^j=af^j$ for $j=1,2$.
Since $a\in L^\infty_{\mathrm{loc}}$, local solvability
of the $\bar\partial$-equation
(see \cite[Chapter I, Sections 4--6]{Vekua})
gives $\sigma\in W^{1,p}_{\mathrm{loc}}$ for every
$p>2$ satisfying $\partial_{\bar z}\sigma=a$.
Consequently,
$$
\partial_{\bar z}(e^{-\sigma}f^j)=0,
\quad j=1,2,
$$
in the weak sense, so
$e^{-\sigma}(f^1,f^2)$ is holomorphic.
A holomorphic
function with zeros accumulating in a connected domain vanishes
identically there. Applying this on overlapping coordinate neighborhoods shows that $\Phi$
either vanishes identically or has only isolated interior zeros. Thus, if
$\Phi\not\equiv0$, at each interior zero $z_0$ we can write
$$
e^{-\sigma}(f^1,f^2)
=(z-z_0)^m h(z),
\quad m\ge1,\quad h(z_0)\ne0,
$$
where $h=(h^1,h^2)$ is holomorphic.
Setting $s=h^1e_1+h^2e_2$ gives
\eqref{eq:Phi-interior-factor}. The same argument applies to $\Psi$.
For an analogous interior factorization in the Weyl-minimal setting,
see~\cite[Section~3.2]{Ream}.

\emph{Boundary.}
Fix a boundary zero of $\Phi$ and choose a
boundary-adapted holomorphic coordinate $\zeta$
centered at that point, with $\partial\Sigma$
represented by the real axis. By
Lemma~\ref{lem:boundary-lagrangian},
$\Phi\in F$ along the boundary.
Choose a smooth complex frame $e_1,e_2$ of $E$ whose boundary
values form a real frame of $F$. Writing $\Phi=\sum_{j=1}^2f^je_j$,
we have $f(\xi)\in\R^2$ for real $\xi$.
Define $A=(A_{ij})$ by
$\bar\partial_Ee_j=\sum_{i=1}^2 A_{ij}e_i\,d\bar\zeta$.
The entries of $A$ are smooth up to the boundary because
$\bar\partial_E$ and the chosen frame are smooth.
Equation~\eqref{eq:scalar-component-equations} then gives
\begin{equation}\label{eq:boundary-system}
\partial_{\bar\zeta}f+Cf=0,
\quad C=A-aI\in L^\infty.
\end{equation}
Extend $f$ and $C$ to the lower half-disk by
$$
f(\zeta)=\overline{f(\bar\zeta)},
\quad
C(\zeta)=\overline{C(\bar\zeta)}
\quad(\operatorname{Im}\zeta<0).
$$
Since $f$ is real-valued on the boundary,
the two pieces have identical Sobolev traces
on the real axis and define a $W^{1,p}$ function
across it. Integration by parts on the two
half-disks shows that the boundary terms cancel.
Thus \eqref{eq:boundary-system} holds weakly
on the full disk.

The complex-linear similarity principle supplies a local
invertible matrix $G\in W^{1,p}$, for some $p>2$, satisfying
$\partial_{\bar\zeta}G+CG=0$; see
\cite[Theorem~2.78 and its proof]{WendlLectures}.
Since $G$ is continuous and $C$ is bounded, local elliptic
regularity gives $G\in W^{1,p}_{\mathrm{loc}}$ for every
finite $p>2$. Thus $h=G^{-1}f$ is weakly holomorphic,
and hence holomorphic. Interior unique continuation excludes
$f\equiv0$ on this neighborhood. Since $f(0)=0$, we may write
$$
h(\zeta)=\zeta^m h_0(\zeta),
\quad m\ge1,\quad h_0(0)\ne0.
$$
Setting $v=Gh_0$, we obtain
$$
f(\zeta)=\zeta^m v(\zeta),
\quad
v\in W^{1,p}_{\mathrm{loc}},
\quad v(0)\ne0.
$$
For real $\xi\ne0$, $v(\xi)=\xi^{-m}f(\xi)\in\R^2$;
continuity gives $v(0)\in\R^2$ as well.
Consequently, $V=\sum_jv^je_j$ is nonvanishing
near the origin and has boundary values in $F$.
This proves \eqref{eq:Phi-boundary-factor}.

The same argument, using a frame adapted to
$\overline F$, applies to $\Psi$.
Hence every zero of a nontrivial $\Phi$ or $\Psi$
is isolated, including at the boundary.
Compactness of $\Sigma$ implies that each zero
set is finite.

The factorizations show that $K_\Phi=\C\Phi$
extends continuously across the zeros of $\Phi$,
with a real boundary line subbundle
$F_\Phi=K_\Phi|_{\partial\Sigma}\cap F$.
The analogous extensions give
$K_\Psi\subset\overline E$ and
$F_\Psi=K_\Psi|_{\partial\Sigma}\cap\overline F$.
\end{proof}

\subsection{Tangency counts and Maslov indices}
\label{subsec:tangency-counts}

Suppose neither $\Phi$ nor $\Psi$ vanishes identically.
For each $x\in\Sigma\setminus\partial\Sigma$,
let $\nu_\Phi(x)$ and $\nu_\Psi(x)$ denote the
vanishing orders of $\Phi$ and $\Psi$, respectively,
with order zero where the corresponding section does not vanish.
Define
\begin{equation}\label{eq:local-RPQ-interior}
\begin{aligned}
R_x&=\min\{\nu_\Phi(x),\nu_\Psi(x)\},\\
P_x&=\nu_\Psi(x)-R_x,\\
Q_x&=\nu_\Phi(x)-R_x.
\end{aligned}
\end{equation}
Since $\mathcal E=E\oplus\overline E$ is an
orthogonal splitting,
\begin{equation}\label{eq:differential-norm}
|u_z|^2=\tfrac14\bigl(|\Phi|^2+|\Psi|^2\bigr).
\end{equation}
Comparison with \eqref{eq:LC-ramification}
shows that $R_x$ is the branching order of $u$
at $x$.

At $q\in\partial\Sigma$, let
$\nu_\Phi^\partial(q)$ and $\nu_\Psi^\partial(q)$
be the reflected orders of $\Phi$ and $\Psi$
in a boundary-adapted holomorphic coordinate.
Define
\begin{equation}\label{eq:local-RPQ-boundary}
\begin{aligned}
R_q^\partial
  &=\min\{\nu_\Phi^\partial(q),
           \nu_\Psi^\partial(q)\},\\
P_q^\partial
  &=\nu_\Psi^\partial(q)-R_q^\partial,\\
Q_q^\partial
  &=\nu_\Phi^\partial(q)-R_q^\partial.
\end{aligned}
\end{equation}
The boundary factorizations and
\eqref{eq:differential-norm} identify
$R_q^\partial$ with the branching order at $q$.

Define the total counts by
\begin{equation}\label{eq:global-RPQ}
\begin{aligned}
R(u)&=\sum_{x\in\Sigma\setminus\partial\Sigma}R_x
 +\tfrac12\sum_{q\in\partial\Sigma}R_q^\partial,\\
P(u)&=\sum_{x\in\Sigma\setminus\partial\Sigma}P_x
 +\tfrac12\sum_{q\in\partial\Sigma}P_q^\partial,\\
Q(u)&=\sum_{x\in\Sigma\setminus\partial\Sigma}Q_x
 +\tfrac12\sum_{q\in\partial\Sigma}Q_q^\partial.
\end{aligned}
\end{equation}
These sums are finite because the zeros are isolated up to the
boundary by Lemma~\ref{lem:relative-factorization} and $\Sigma$ is compact.
Interior orders are counted in full, while boundary orders are counted with
weight $1/2$. Thus the total counts may be
half-integers.

If $Z(\Phi)$ and $Z(\Psi)$ denote the corresponding
weighted zero counts before removing ramification,
then
\begin{equation}\label{eq:raw-divisors}
Z(\Phi)=R(u)+Q(u),\quad
Z(\Psi)=R(u)+P(u).
\end{equation}

At a branch point, the oriented tangent plane
extends continuously. Indeed, dividing $u_z$
by its common local factor yields a continuous
section $V$ with $V\ne0$ at the branch point.
Conformality implies $(V,V)=0$, so
$\operatorname{Re}V$ and $-\operatorname{Im}V$
span an oriented two-plane. Multiplication by
a nonzero complex scalar preserves this plane
and its orientation. Consequently, the
K\"ahler angle extends continuously, with
\begin{equation}\label{eq:angle-leading}
\cos\alpha
=\frac{|V^{1,0}|^2-|V^{0,1}|^2}
       {|V^{1,0}|^2+|V^{0,1}|^2}.
\end{equation}
At a branch point,
\begin{align*}
P_x>0
&\iff V^{0,1}(x)=0
\iff\cos\alpha(x)=1,\\
Q_x>0
&\iff V^{1,0}(x)=0
\iff\cos\alpha(x)=-1.
\end{align*}
The same equivalences hold at boundary branch
points, using the corresponding boundary orders.
Thus $P(u)$ and $Q(u)$ count complex and
anticomplex tangencies after removal of the
common ramification; compare
Chen--Tian~\cite[pp.~902--903]{ChenTian}.

Under a holomorphic change $z=z(\zeta)$, both $\Phi$ and $\Psi$
are multiplied by the nonvanishing factor $z'(\zeta)$.
Their vanishing orders, and hence all the counts above, are
coordinate independent.

If $u$ is $J$-holomorphic, then $\Psi\equiv0$.
In this case, set $R(u)=Z(\Phi)$ and $Q(u)=0$,
without assigning a finite value to $P(u)$.
If $u$ is $J$-antiholomorphic, set
$R(u)=Z(\Psi)$ and $P(u)=0$;
$Q(u)$ is not defined in this case.

\begin{lemma}\label{lem:line-zero-count}
For a bordered surface $\Sigma$ as in Lemma~\ref{lem:relative-factorization},
if $\Phi\not\equiv0$, then
\begin{equation}\label{eq:Phi-line-Maslov}
\mu(K_\Phi,F_\Phi)=2\bigl(Z(\Phi)+\chi(\Sigma)\bigr).
\end{equation}
If $\Psi\not\equiv0$, then
\begin{equation}\label{eq:Psi-line-Maslov}
\mu(K_\Psi,F_\Psi)=2\bigl(Z(\Psi)+\chi(\Sigma)\bigr).
\end{equation}
\end{lemma}

\begin{proof}
Recall that $K_\Sigma=(T^{1,0}\Sigma)^*$ is the canonical line bundle.
Let $G_\Sigma\subset K_\Sigma|_{\partial\Sigma}$ be the real line
subbundle with fibers
$$
G_{\Sigma,p}=\{\beta\in K_{\Sigma,p}:\beta(v)\in\R
\text{ for every }v\in T_p\partial\Sigma\},\quad p\in\partial\Sigma.
$$
In a boundary-adapted holomorphic coordinate $z=s+it$ with
$\partial\Sigma$ represented by $t=0$, a one-form $\beta=c\,dz$
satisfies $\beta(\partial_s)=c$. Hence $G_\Sigma$ is locally
spanned over $\R$ by $dz$ along the boundary.
Let $\widehat\Sigma$ denote the closed oriented double of $\Sigma$,
obtained by gluing two copies along corresponding boundary points,
with the second copy given the opposite orientation.
The tangent bundle pair $(T\Sigma,T\partial\Sigma)$, with
$T\Sigma$ viewed as a complex line bundle via $j$, has Maslov index
$2\chi(\Sigma)$: its double is the tangent bundle of
$\widehat\Sigma$, whose Chern number is
$\chi(\widehat\Sigma)=2\chi(\Sigma)$.
The pair $(K_\Sigma,G_\Sigma)$ is the dual of
$(T\Sigma,T\partial\Sigma)$. Its double is therefore the dual
of the tangent bundle of $\widehat\Sigma$.
Since a complex line bundle and its dual have opposite first
Chern classes, we obtain
$$
\mu(K_\Sigma,G_\Sigma)
=-\mu(T\Sigma,T\partial\Sigma)
=-2\chi(\Sigma).
$$

Suppose $\Phi\not\equiv0$. By
Lemma~\ref{lem:boundary-lagrangian}, the global section
$\boldsymbol\Phi=\Phi\,dz$ of $K_\Sigma\otimes K_\Phi$
takes boundary values in $G_\Sigma\otimes F_\Phi$.
Its local zero indices are positive by
Lemma~\ref{lem:relative-factorization}. On doubling this
line bundle pair and section, each interior zero of order
$m$ appears twice, while each boundary zero of reflected
order $m$ appears once. The latter assertion follows in
the boundary frame $V$ from $\Phi=\zeta^mV$:
reflection extends the scalar coefficient as $\zeta^m$.
Thus the doubled section has total zero index $2Z(\Phi)$.
This topological zero count applies to the continuous
extended line bundles; it does not require a smooth or
holomorphic extension across the boundary. See
\cite[Appendix, Propositions~A.1--A.2]{Wendl}.

The Chern number of the doubled line bundle equals the
Maslov index of the original pair. By tensor-product
additivity for line bundle pairs,
$$
\begin{aligned}
2Z(\Phi)
&=\mu(K_\Sigma\otimes K_\Phi,G_\Sigma\otimes F_\Phi)\\
&=-2\chi(\Sigma)+\mu(K_\Phi,F_\Phi).
\end{aligned}
$$
This proves \eqref{eq:Phi-line-Maslov}. Applying the same
argument to $\boldsymbol\Psi$ proves
\eqref{eq:Psi-line-Maslov}.
\end{proof}

\section{Lichnerowicz homotopy identity under a Lagrangian constraint}
\label{Lichnerowicz}
We extend Lichnerowicz's partial-energy homotopy identity
\cite{Lichnerowicz} to compact almost K\"ahler domains with boundary
mapped into a Lagrangian submanifold. The boundary constraint
eliminates the boundary term and gives energy minimization in
relative homotopy classes in arbitrary dimension.

Let $(N^{2n},g_N,J_N,\omega_N)$ be a compact connected almost Hermitian
manifold with smooth boundary, and let $(M,g_M,J_M,\omega_M)$ be almost
Hermitian.  For a smooth
map $\varphi:N\to M$, decompose the complexified differential as
\begin{equation}\label{eq:Lichnerowicz-partial-differentials}
\begin{aligned}
\partial\varphi &:T^{1,0}N\longrightarrow T^{1,0}M,&
\bar\partial\varphi &:T^{0,1}N\longrightarrow T^{1,0}M,\\
\partial\bar\varphi &:T^{1,0}N\longrightarrow T^{0,1}M,&
\bar\partial\bar\varphi &:T^{0,1}N\longrightarrow T^{0,1}M.
\end{aligned}
\end{equation}
With the induced Hermitian norms, set 
\begin{equation}\label{eq:Lichnerowicz-partial-energies}
e'(\varphi)=|\partial\varphi|^2,\quad
e''(\varphi)=|\bar\partial\varphi|^2,\quad
e(\varphi)=e'(\varphi)+e''(\varphi)=\tfrac12|d\varphi|^2, 
\end{equation}
and
$$
E'(\varphi)=\int_N e'(\varphi)\,dv_{g_N},\quad 
E''(\varphi)=\int_N e''(\varphi)\,dv_{g_N},\quad E(\varphi)=E'(\varphi)+E''(\varphi).
$$ 

\begin{lemma}
\label{lem:Lichnerowicz-homotopy}
Assume that $N$ and $M$ are almost K\"ahler manifolds and that
$L\subset M$ is Lagrangian.  For every smooth relative homotopy
$\varphi_t:(N,\partial N)\to(M,L)$, the quantity
$E'(\varphi_t)-E''(\varphi_t)$ is constant.
If $\omega_M=d_M\lambda$ and $\lambda|_L=d_L f$, then
$$
E'(\varphi_t)=E''(\varphi_t)=\tfrac12E(\varphi_t).
$$
\end{lemma}

\begin{proof}
Let $*_N$ be the Hodge star of $g_N$ and set
$$
\theta=*_N\omega_N=\frac{\omega_N^{n-1}}{(n-1)!},
\quad d_N\theta=0.
$$
Lichnerowicz's pointwise identity~\cite{Lichnerowicz} gives
$$
\bigl(e'(\varphi_t)-e''(\varphi_t)\bigr)\,dv_{g_N}
=\varphi_t^*\omega_M\wedge\theta.
$$
Define the variation one-form $\beta_t$ on $N$ by
$$
\beta_t(X)=\omega_M(\dot\varphi_t,d\varphi_tX).
$$
Since $d_M\omega_M=0$, the homotopy formula gives
$\frac{d}{dt}(\varphi_t^*\omega_M)=d_N\beta_t$;
see~\cite[Homotopy Lemma (8.8)]{E-L}. Writing
$\iota:\partial N\hookrightarrow N$ for the inclusion, we obtain
\begin{equation}\label{t derivative}
\begin{aligned}
\frac{d}{dt}\bigl(E'(\varphi_t)-E''(\varphi_t)\bigr)
&=\int_N d_N\beta_t\wedge\theta\\
&=\int_N d_N(\beta_t\wedge\theta)
=\int_{\partial N}\iota^*(\beta_t\wedge\theta).
\end{aligned}
\end{equation}
For $x\in\partial N$ and $X\in T_x\partial N$, both
$\dot\varphi_t(x)$ and $d\varphi_tX$ lie in $T_{\varphi_t(x)}L$.
The Lagrangian condition therefore gives $\iota^*\beta_t=0$,
proving homotopy invariance.

Now suppose $\omega_M=d_M\lambda$ and $\lambda|_L=d_Lf$.
Put $f_t=f\circ\varphi_t|_{\partial N}$. Then
$\iota^*\varphi_t^*\lambda=d_{\partial N}f_t$, so Stokes' theorem
and $d_N\theta=0$ give
$$
\begin{aligned}
E'(\varphi_t)-E''(\varphi_t)
&=\int_N d_N(\varphi_t^*\lambda)\wedge\theta\\
&=\int_{\partial N}d_{\partial N}f_t\wedge\iota^*\theta\\
&=\int_{\partial N}d_{\partial N}(f_t\,\iota^*\theta)=0.
\end{aligned}
$$
Since $E=E'+E''$, it follows that $E'=E''=E/2$.
\end{proof}

\begin{lemma}
\label{lem:Lichnerowicz-holomorphic-minimizer}
Under the hypotheses of Lemma~\ref{lem:Lichnerowicz-homotopy}, a
$J$-holomorphic or $J$-antiholomorphic map
$u:(N,\partial N)\to(M,L)$ minimizes energy in its relative homotopy
class. Every energy minimizer in that class is holomorphic,
respectively antiholomorphic. If $\omega_M=d_M\lambda$ and
$\lambda|_L$ is exact, every such
holomorphic or antiholomorphic map is constant.
\end{lemma}

\begin{proof}
Suppose $u$ is $J$-holomorphic, so $E''(u)=0$ and $E(u)=E'(u)$.
For every smooth map $v$ relatively homotopic to $u$,
Lemma~\ref{lem:Lichnerowicz-homotopy} gives
$$
E(v)=\bigl(E'(v)-E''(v)\bigr)+2E''(v)
=E(u)+2E''(v)\ge E(u).
$$
Equality holds precisely when $E''(v)=0$, or equivalently
$\bar\partial v=0$, since $v$ is smooth. Thus every energy minimizer
in the class is $J$-holomorphic. If $u$ is $J$-antiholomorphic,
interchanging $E'$ and $E''$ gives $E(v)=E(u)+2E'(v)$, with equality
precisely when $v$ is $J$-antiholomorphic.

In the exact case, the preceding lemma gives $E'=E''$.  One of these
partial energies vanishes for a holomorphic or antiholomorphic map,
so $E=0$.  Therefore $du=0$, and connectedness of $N$ gives constancy.
\end{proof}

In dimension two, Lemma~\ref{lem:Lichnerowicz-homotopy} and
Wirtinger's inequality give the following consequence.

\begin{corollary}\label{cor:Lichnerowicz-surface}
Let $(\Sigma,j)$ be a compact connected Riemann surface with smooth
boundary, let $(M,g,J,\omega)$ be almost K\"ahler, and let $L\subset M$
be Lagrangian.  A $J$-holomorphic or $J$-antiholomorphic map
$u:(\Sigma,\partial\Sigma)\to(M,L)$ minimizes area in its
relative homotopy class and satisfies the free-boundary condition.
\end{corollary}

\begin{proof}
By Wirtinger's inequality (see
\cite[Section~6, Example~I]{HarveyLawson}),
for any real tangent vectors $X,Y\in TM$,
$$
|\omega(X,Y)|\le |X\wedge Y|_g
=\sqrt{|X|_g^2|Y|_g^2-g(X,Y)^2}.
$$
For linearly independent $X,Y$, equality holds precisely
when the plane spanned by $X,Y$ is $J$-invariant.
In particular, for a $J$-holomorphic or
$J$-antiholomorphic map $u$,
$$
\Area(u)=\left|\int_\Sigma u^*\omega\right|.
$$
This equality remains valid at branch points.

Fix an arbitrary smooth $j$-compatible background metric $h$ on
all of $\Sigma$. Every $j$-holomorphic coordinate $z=x+iy$ is
isothermal for $h$, so $h=e^{2\rho}(dx^2+dy^2)$.
Since $j$ is integrable and the fundamental two-form
$\omega_\Sigma$ is closed in real dimension two, $(\Sigma,h,j)$
is K\"ahler. Use $h$ to define the partial energies.
Let $v:(\Sigma,\partial\Sigma)\to(M,L)$ be any smooth map
relatively homotopic to $u$, and choose a smooth relative homotopy
$u_t$ with $u_0=u$ and $u_1=v$. For this domain, $n=1$ and
$\theta=*_\Sigma\omega_\Sigma=1$ in the proof of
Lemma~\ref{lem:Lichnerowicz-homotopy}. The partial-energy identity
therefore reads $E'(u_t)-E''(u_t)=\int_\Sigma u_t^*\omega$.
Lemma~\ref{lem:Lichnerowicz-homotopy} states that the left-hand side
is independent of $t$. Evaluating at the endpoints gives
$$
\int_\Sigma v^*\omega
=E'(v)-E''(v)
=E'(u)-E''(u)
=\int_\Sigma u^*\omega.
$$
Applying Wirtinger's inequality pointwise,
we obtain
$$
\Area(v)
\ge\left|\int_\Sigma v^*\omega\right|
=\left|\int_\Sigma u^*\omega\right|
=\Area(u).
$$
Hence $u$ minimizes area in its relative homotopy class.
Consequently, the first variation of area vanishes for every smooth
variation $u_t$ satisfying $u_0=u$ and
$u_t(\partial\Sigma)\subset L$ for all sufficiently small $|t|$.
The variational field is tangent to $L$ along the boundary, and
its boundary values may be chosen arbitrarily within $TL$ near
any regular boundary point. The first variation formula therefore
gives the free-boundary condition
$du_p(\nu_p)\perp T_{u(p)}L$ at every regular boundary point $p$,
where $\nu$ is the outward unit normal to $\partial\Sigma$
for $u^*g$.
\end{proof}

\section{Complex-point formula and rigidity}\label{sec:rigidity}

\subsection{Complex-point formula and holomorphicity criterion}
\label{sec:canonical-line}

\begin{theorem}\label{thm:main}\label{thm:AK-main}
Let $(M^4,g,J,\omega)$ be almost K\"ahler, let $L^2\subset M$ be a
smooth oriented embedded Lagrangian surface, and let
$u:(\Sigma,\partial\Sigma)\to(M,L)$ be a smooth nonconstant branched
minimal immersion satisfying the free-boundary condition. Then the following hold.
\begin{enumerate}[label=(\roman*)]
\item If $u$ is neither $J$-holomorphic nor $J$-antiholomorphic, then
\begin{equation}\label{eq:AK-two-sided-Maslov}
\mu_L(u)=2\bigl(Q(u)-P(u)\bigr).
\end{equation}
\item If $u$ is not $J$-antiholomorphic and
\begin{equation}\label{eq:AK-zero-count-criterion}
\mu_L(u)>2Q(u),
\end{equation}
then $u$ is $J$-holomorphic and minimizes area in its relative
homotopy class. In particular, the conclusion holds when $u$
has no anticomplex points and $\mu_L(u)>0$.
\end{enumerate}
\end{theorem}

\begin{proof}
\emph{(i)} Neither $\Phi$ nor $\Psi$ vanishes identically.
By Lemma~\ref{lem:relative-factorization}, their complex
line subbundles $K_\Phi\subset E$ and
$K_\Psi\subset\overline E$ extend continuously across all
zeros, including boundary zeros. Their real boundary lines
are $F_\Phi=K_\Phi\cap F$ and
$F_\Psi=K_\Psi\cap\overline F$. 
Conformality and the isotropy of $E$ and $\overline E$
for the complex-bilinear pairing give
$$
0=(u_z,u_z)=\tfrac12(\Phi,\Psi)
 =\tfrac12\langle\Phi,\overline\Psi\rangle.
$$
Where both components are nonzero, $K_\Phi$ and
$\overline{K_\Psi}$ are therefore orthogonal complex lines
in $E$. Their zero sets are finite, so this orthogonality
extends by continuity to every point of $\Sigma$.
Since $\rank_\C E=2$, it follows that
\begin{equation}\label{eq:E-differential-splitting}
E=K_\Phi\oplus\overline{K_\Psi}.
\end{equation}
Along the boundary, $F_\Phi$ and $\overline{F_\Psi}$ are
real lines in $F$ lying in these orthogonal summands.
They are independent, including at boundary zeros, and span $F$
because $\rank_\R F=2$. Consequently,
\begin{equation}\label{eq:EF-splitting}
(E,F)=(K_\Phi,F_\Phi)
\oplus(\overline{K_\Psi},\overline{F_\Psi}).
\end{equation}

Conjugation reverses the Maslov index. To check the sign,
choose a complex trivialization of $K_\Psi$, which exists
because $\Sigma$ has nonempty boundary. On each boundary
component $C$, with its induced orientation, write
$F_\Psi=e^{i\vartheta_C}\R$, where
$\vartheta_C:C\to\R/\pi\Z$. Then
$$
\mu(K_\Psi,F_\Psi)
=\sum_C\wind_C(e^{2i\vartheta_C}).
$$
In the conjugate trivialization the doubled phase is
$e^{-2i\vartheta_C}$. Thus
$$
\mu(\overline{K_\Psi},\overline{F_\Psi})
=-\mu(K_\Psi,F_\Psi).
$$
The phase is taken modulo $\pi$, so this calculation also applies
when $F_\Psi$ is not orientable.

Maslov additivity in \eqref{eq:EF-splitting}, together with
Lemma~\ref{lem:line-zero-count} and \eqref{eq:raw-divisors}, gives
$$
\begin{aligned}
\mu_L(u)
&=\mu(K_\Phi,F_\Phi)-\mu(K_\Psi,F_\Psi)\\
&=2\bigl(R(u)+Q(u)+\chi(\Sigma)\bigr)
 -2\bigl(R(u)+P(u)+\chi(\Sigma)\bigr)\\
&=2\bigl(Q(u)-P(u)\bigr).
\end{aligned}
$$
In particular, both the ramification and
Euler-characteristic terms cancel. This proves (i).

\emph{(ii)} If $u$ were not $J$-holomorphic, the hypothesis
that it is not $J$-antiholomorphic would allow us to apply (i).
Since $P(u)\ge0$, we would obtain
$$
\mu_L(u)=2\bigl(Q(u)-P(u)\bigr)\le2Q(u),
$$
contrary to \eqref{eq:AK-zero-count-criterion}.
Hence $u$ is $J$-holomorphic. Area minimization
then follows from Corollary~\ref{cor:Lichnerowicz-surface}.
If $u$ has no anticomplex points, $Q(u)=0$, which gives
the final assertion.
\end{proof}

\subsection{K\"ahler--Einstein rigidity}\label{subsec:KE-angle}
We first combine the relative complex-point formula with
\eqref{eq:maslov-integral} to prove rigidity for positive Einstein
constant, and then we use the K\"ahler-angle equation for the hyperk\"ahler case.

\begin{theorem}\label{thm:KE-positive}
Let $(M^4,g,J,\omega)$ be K\"ahler--Einstein with
$\Ric_M=\kappa g$, $\kappa>0$. Let $L^2\subset M$ be an oriented
minimal Lagrangian surface.  Let
$u:(\Sigma,\partial\Sigma)\to(M,L)$ be a smooth nonconstant branched
minimal immersion satisfying the free-boundary condition, with
K\"ahler angle $\alpha$.
Then, 
\begin{enumerate}[label=(\alph*)]
\item For each $\varepsilon\in\{+1,-1\}$, if
$\varepsilon\cos\alpha\ge0$ on $\Sigma$, then either $u$ is 
$\varepsilon J$-holomorphic and minimizes area in its relative
homotopy class, or $u$ is Lagrangian. In particular, if $u$ has no
Lagrangian points, then it is $J$-holomorphic or
$J$-antiholomorphic and minimizes area in its relative homotopy class.
\item If $u$ is totally real, then $\omega(u)=\mu_L(u)=0$.
\item For every $A\in H_2(M,L;\Z)$,
\begin{equation}\label{eq:KE-Maslov}
\mu_L(A)=\frac{\kappa}{\pi}\omega(A),
\quad \mu_L(A)\in2\Z.
\end{equation}
\end{enumerate}
\end{theorem}

\begin{proof}
We first prove (c). Let $\Sigma'$ be any compact oriented surface,
possibly with boundary, and let
$v:(\Sigma',\partial\Sigma')\to(M,L)$ be smooth.
Let $\rho_M(X,Y)=\Ric_M(JX,Y)$ be the Ricci form.
By \cite{CG} and \cite[Corollary~4.2]{Pacini},
\begin{equation}\label{eq:maslov-integral}
\mu_L(v)=\frac1\pi\left(
\int_{\Sigma'}v^*\rho_M
-\int_{\partial\Sigma'}(v|_{\partial\Sigma'})^*\omega(H_L,\cdot)
\right),
\end{equation}
where $H_L=\tr_{TL}B^L$. Since $H_L=0$ and
$\rho_M=\kappa\omega$, this gives
$\mu_L(v)=(\kappa/\pi)\omega(v)$.
Both sides depend only on the relative homology class
$A=v_*[\Sigma',\partial\Sigma']$.  By additivity and the representation
of relative two-dimensional classes by sums of mapped oriented surfaces,
the relation holds for every $A\in H_2(M,L;\Z)$.
Since $\mathfrak m_L$ is twice an integral relative Chern class,
$\mu_L(A)\in2\Z$, giving
the asserted symplectic-area periods.

To prove (a), suppose first that $\cos\alpha\ge0$ on $\Sigma$.
There are no anticomplex points, so $Q(u)=0$.
If $u$ is $J$-holomorphic, area minimization follows from
Corollary~\ref{cor:Lichnerowicz-surface}.
Otherwise, $u$ is neither $J$-holomorphic nor $J$-antiholomorphic,
and Theorem~\ref{thm:AK-main}(i) and \eqref{eq:KE-Maslov} give
$$
0\le \frac{\kappa}{\pi}\int_\Sigma\cos\alpha\,dA_u
=\mu_L(u)=2\bigl(Q(u)-P(u)\bigr)=-2P(u)\le0.
$$
Since $\cos\alpha$ is continuous and nonnegative, it vanishes
identically. Thus $u$ is Lagrangian.
The case $\cos\alpha\le0$ follows by replacing $J$ with $-J$.

For (b), total reality gives $P(u)=Q(u)=0$. The map is
neither holomorphic nor antiholomorphic, so
Theorem~\ref{thm:AK-main}(i) gives
$\mu_L(u)=0$.  Equation~\eqref{eq:KE-Maslov} then implies 
$\omega(u)=0$.
\end{proof}

The argument is the free-boundary counterpart of Webster's
closed-surface complex-point formula~\cite{Webster,WebsterChern}.
For a closed oriented minimal surface $v:\Sigma_0\to M$ that is
neither holomorphic nor antiholomorphic, that formula reads
$$
Q(v)-P(v)=\langle c_1(TM),v_*[\Sigma_0]\rangle
=\frac{\kappa}{2\pi}\int_{\Sigma_0}\cos\alpha\,dA_v.
$$
Theorem~\ref{thm:AK-main}(i) and \eqref{eq:KE-Maslov} give the
same relation in the present boundary problem, with the Chern number
replaced by one half of the relative Maslov index and boundary
tangencies counted by half their reflected orders. In both settings,
positivity of $\kappa$ gives rigidity when $\cos\alpha$ has one sign.

When $\kappa=0$, equation~\eqref{eq:maslov-integral} gives $\mu_L(u)=0$.
Theorem~\ref{thm:AK-main}(i) then gives $P(u)=Q(u)$ for a surface
that is not $\pm$ holomorphic, but does not
control the K\"ahler angle. The next theorem uses the angle equation, following
Wolfson's closed-surface arguments~\cite[Theorems~2.1 and~2.2]{Wolfson}.

\begin{theorem}\label{thm:KE-boundary-angle}\label{thm:KE-hyperkahler}
Let $(M^4,g,J,\omega)$ be hyperk\"ahler, and let
$L^2\subset M$ be an oriented minimal Lagrangian surface.
Let $u:(\Sigma,\partial\Sigma)\to(M,L)$ be a smooth nonconstant
branched minimal immersion satisfying the free-boundary condition.
If $u$ has no Lagrangian points on $\Sigma$, then $u$ is
holomorphic for one fixed parallel complex structure in the
hyperk\"ahler $S^2$-family.
\end{theorem}

\begin{proof}
Let $\mathcal B$ be the finite branch set and fix an arbitrary smooth
$j$-compatible background metric $h$ on all of $\Sigma$. Locally, $h=e^{2\rho}(dx^2+dy^2)$
in a $j$-holomorphic (hence isothermal) coordinate $z=x+iy$.
Write $u^*g=e^{2\varphi}h$ on the regular set. The factor $e^{2\varphi}$ extends continuously to $\Sigma$,
vanishing precisely at $\mathcal B$; geometric quantities for the
induced metric are understood on $\Sigma\setminus\mathcal B$.

\emph{Neumann boundary condition.}
Let $\nu_u$ be the outward unit normal to $\partial\Sigma$ for $u^*g$
and let $\eta=du(\nu_u)$ be the corresponding unit conormal of
the immersed surface. Choose the unit boundary tangent $T$
so that $(T,\eta)$ is positively oriented. In the calculation below,
tangent vectors on the regular domain are identified with their
images under $du$.
Choose a unit vector $W\in TL$, $W\perp T$, with
$$
\eta=\cos\alpha\,JT+\sin\alpha\,JW,
\quad
\xi=-\sin\alpha\,JT+\cos\alpha\,JW.
$$
Then $(W,\xi)$ is an orthonormal normal frame of the surface.
Extend $T,\eta$ locally with
$\nabla^\Sigma_\eta T=\nabla^\Sigma_\eta\eta=0$ at the point.
Minimality of $u$ gives
$B^\Sigma(\eta,\eta)=-B^\Sigma(T,T)$; consequently,
$$
\begin{aligned}
\partial_\eta\cos\alpha
&=g(JB^\Sigma(\eta,T),\eta)-g(JT,B^\Sigma(T,T))\\
&=\sin\alpha\,(a+b),
\end{aligned}
$$
where
$a=g(B^\Sigma(T,\eta),W)$ and
$b=g(B^\Sigma(T,T),\xi)$. The free-boundary identities give
$$
a=-g(\eta,B^L(T,W)),
\quad b=g(B^L(T,T),\xi).
$$
Since $L$ is Lagrangian, the cubic form
$C(X,Y,Z)=g(B^L(X,Y),JZ)$ is totally symmetric for
$X,Y,Z\in TL$.
Substitution therefore yields
$$
a+b=-\sin\alpha\,g(\tr_{TL}B^L,JT)=0,
$$
because $L$ is minimal. Thus
\begin{equation}\label{eq:angle-boundary-Neumann}
\partial_{\nu_u}\cos\alpha=\partial_\eta\cos\alpha=0
\end{equation}
on each regular boundary arc away from branch points and complex
and anticomplex points.

\emph{Angle equation and exceptional points.}
For a K\"ahler--Einstein metric with $\Ric_M=\kappa g$,
 the known interior K\"ahler-angle equation is
valid away from branch points and complex and anticomplex points:
\begin{equation}\label{eq:KE-angle-interior}
\Delta_\Sigma\cos\alpha
=-\cos\alpha\bigl(2|d\alpha|_\Sigma^2
+\kappa\sin^2\alpha\bigr),
\end{equation}
where $\Delta_\Sigma$ and
$|d\alpha|_\Sigma$ are computed with the induced metric $u^*g$; see~\cite[pp.~904--905]{ChenTian} and~\cite{Wolfson}. 

Suppose $u$ is not $\pm J$-holomorphic.
The complex and anticomplex point sets are finite by
Lemma~\ref{lem:relative-factorization}. Let $\mathcal S$ be their
union with $\mathcal B$. Near a branch point of order $r$, the
common branching factor can be removed locally:
$$
\Phi=\zeta^r\widehat\Phi,
\quad
\Psi=\zeta^r\widehat\Psi.
$$
The factors $\widehat\Phi,\widehat\Psi$ do not vanish simultaneously.
They are smooth in the interior and $W^{1,p}$ up to the boundary for
every finite $p>2$, after the boundary reflection.  Thus
\begin{equation}\label{eq:angle-branch-factor}
\cos\alpha
=\frac{|\widehat\Phi|^2-|\widehat\Psi|^2}
       {|\widehat\Phi|^2+|\widehat\Psi|^2}
\end{equation}
extends across each branch point.  Where $\widehat\Phi\ne0$,
$$
\alpha=2\arctan\frac{|\widehat\Psi|}{|\widehat\Phi|}.
$$
Similarly, where $\widehat\Psi\ne0$, we have
$\pi-\alpha=2\arctan(|\widehat\Phi|/|\widehat\Psi|)$.
The local factorizations at regular complex and
anticomplex points, including their reflected boundary versions,
give the same regularity.  Therefore $\cos\alpha$ and $\alpha$ belong
to $W^{1,p}$ locally across $\mathcal S$ for every $p>2$.

Whenever $\cos\alpha$ has a fixed sign, set
$f=\varepsilon\cos\alpha\ge0$ with
$\varepsilon\in\{1,-1\}$ chosen accordingly; zeros of
$\cos\alpha$ are allowed. Since $u^*g=e^{2\varphi}h$ on the regular
set, the two-dimensional conformal transformation laws give
$\Delta_{u^*g}=e^{-2\varphi}\Delta_h$ and
$|d\alpha|_{u^*g}^2=e^{-2\varphi}|d\alpha|_h^2$.
Thus the K\"ahler-angle equation becomes
\begin{equation}\label{eq:angle-parameter}
\Delta_h f=-2f|d\alpha|_h^2
-\kappa f\sin^2\alpha\,e^{2\varphi}
\quad\text{on }\Sigma\setminus\mathcal S.
\end{equation}
Let $\nu_h$ be the outward $h$-unit normal to $\partial\Sigma$.
Since $u^*g=e^{2\varphi}h$, the unit normals satisfy
$\nu_u=e^{-\varphi}\nu_h$. Thus
\eqref{eq:angle-boundary-Neumann} gives
$$
\partial_{\nu_h}f
=e^{\varphi}\partial_{\nu_u}f
=\varepsilon e^{\varphi}\partial_\eta\cos\alpha=0
$$
on the regular boundary arcs under consideration. Hence the
homogeneous Neumann condition holds for $h$ as well as for $u^*g$.

The finite set $\mathcal S$ can be removed from integration by parts.
Indeed, in every interior or boundary coordinate chart the local
factorizations give $\alpha,f\in W^{1,p}$ for any finite $p>2$.
Choose $q=p/(p-1)<2$ and smooth cutoffs $\chi_\delta$ on $\Sigma$
that vanish near $\mathcal S$, equal one away from shrinking
neighborhoods of $\mathcal S$, and satisfy
$$
\|d\chi_\delta\|_{L^q(\Sigma,h)}\longrightarrow0
\quad\text{as }\delta\downarrow0.
$$
Integrate \eqref{eq:angle-parameter} against
$\chi_\delta$ with respect to $dA_h$.  The boundary
term vanishes by the Neumann condition and H\"older's inequality gives
$$
\left|\int_\Sigma\langle df,d\chi_\delta\rangle_h\,dA_h\right|
\le \|df\|_{L^p(h)}\|d\chi_\delta\|_{L^q(h)}
\longrightarrow0.
$$
The remaining integrands are in $L^1$: in particular
$|d\alpha|_h^2\in L^{p/2}$.
Passing to the limit yields the coordinate-independent identity
\begin{equation}\label{eq:angle-integral}
\begin{aligned}
0=\int_\Sigma f\bigl(2|d\alpha|_h^2
+\kappa\sin^2\alpha\,e^{2\varphi}\bigr)\,dA_h
 =\int_\Sigma f\bigl(2|d\alpha|_{u^*g}^2
+\kappa\sin^2\alpha\bigr)\,dA_u.
\end{aligned}
\end{equation}
The second integral is over the regular set. The cutoff estimate
shows that the exceptional points contribute no defect term.

\emph{Hyperk\"ahler rotation.}
Let $\kappa=0$ and assume that $u$ has no Lagrangian points.
If $u$ is $J$-holomorphic or $J$-antiholomorphic, take $J_0=J$ or
$J_0=-J$. Otherwise $f=|\cos\alpha|$ is continuous and strictly
positive on $\Sigma$, and \eqref{eq:angle-integral} gives
$d\alpha=0$ on $\Sigma\setminus\mathcal S$.
This set is connected, so $\alpha$ is constant there and hence,
by continuity, on all of $\Sigma$. Since $u$ is not
$\pm J$-holomorphic, the constant lies in $(0,\pi)$.

We conclude by the hyperk\"ahler rotation argument; compare
\cite[Theorem~2.2]{Wolfson} and the proof of
\cite[Theorem~5.3]{MW}. Choose a parallel quaternionic triple
$J_1=J,J_2,J_3$ and equip
its real span with the normalized inner product
$\langle A,B\rangle_{\mathrm{HK}}=-\tfrac14\tr_{\R}(AB)$,
so that $\langle J_i,J_j\rangle_{\mathrm{HK}}=\delta_{ij}$.
For a positively oriented orthonormal tangent frame $(e_1,e_2)$,
the self-dual Gauss map in this parallel trivialization is
$$
G=\sum_{i=1}^3 g(J_i e_1,e_2)J_i.
$$
Since $J_1e_1,J_2e_1,J_3e_1$ form an orthonormal basis of
$e_1^\perp$, we have $|G|_{\mathrm{HK}}=1$ and $Ge_1=e_2$.
In particular (cf.~\cite[equation~(4)]{Arezzo}),
$$
\langle G,J\rangle_{\mathrm{HK}}=g(Je_1,e_2)=\cos\alpha.
$$
For a conformal minimal immersion, $G$ is holomorphic or
antiholomorphic; see~\cite[Section~2]{Arezzo}.
Since $\alpha$ is constant, the image of $G$ lies in a latitude
circle of $S^2$. The open-mapping theorem forces $G$ to be
constant, say $G\equiv J_0$. Thus $J_0$ is one fixed parallel
complex structure and $du\circ j=J_0\circ du$ on the regular surface.
The identity extends across branch points by continuity.
\end{proof}

\begin{proof}[Alternative proof of Theorem~\ref{thm:KE-positive}(a)]
Suppose $\kappa>0$. The holomorphic and antiholomorphic cases
follow from Corollary~\ref{cor:Lichnerowicz-surface}.
Otherwise, under the sign hypothesis,
$f=|\cos\alpha|$ satisfies \eqref{eq:angle-integral}, including
at Lagrangian points. On the open set where $f>0$, the
nonnegative integrand and $\kappa>0$ force $\sin\alpha=0$,
hence $f=1$.  By continuity, $f$ takes values only
in $\{0,1\}$ on the connected bordered surface, so it is constant.
The case $f\equiv1$ would make $u$ holomorphic or antiholomorphic,
contrary to the case under consideration. Thus $f\equiv0$ and $u$
is Lagrangian.
\end{proof}

\begin{corollary}
\label{cor:KE-period-tangency}
Let $(M^4,g,J,\omega)$ be K\"ahler--Einstein with
$\Ric_M=\kappa g$, $\kappa>0$, and let $L^2\subset M$ be an oriented
minimal Lagrangian surface.  If
$u:(\Sigma,\partial\Sigma)\to(M,L)$ is a branched minimal immersion
satisfying the free-boundary condition and neither $J$-holomorphic
nor $J$-antiholomorphic, then
\begin{equation}\label{eq:KE-period-tangency}
\omega(u)=\frac{2\pi}{\kappa}\bigl(Q(u)-P(u)\bigr).
\end{equation}
\end{corollary}
\begin{proof}
Combine \eqref{eq:KE-Maslov} with Theorem~\ref{thm:AK-main}(i).
\end{proof}
\begin{remark}
In particular, zero symplectic area is equivalent to $P(u)=Q(u)$.
Positive symplectic area gives $Q(u)>P(u)$, and negative symplectic
area gives $P(u)>Q(u)$.
\end{remark}

\section{Examples}\label{sec:examples-sharpness}
The following examples illustrate the sharpness and limitations of our main results.

\subsection{An unstable minimal disk with Lagrangian free boundary}
\label{subsec:CP2-nonminimizing}
Equip $\mathbb{CP}^2$ with the Fubini--Study metric normalized so
that a projective line has area $\pi$. Let
$$
L=\big\{[Z_0:Z_1:Z_2]:|Z_0|=|Z_1|=|Z_2|\big\}
$$
be the Clifford torus, which is minimal and Lagrangian.  Consider $u_0:\D\to\mathbb{CP}^2$ defined by 
$$
u_0(z)=[1:z:\bar z],\quad |z|\le1.
$$
Its image lies in the fixed set of the antiholomorphic isometric
involution
$$
[Z_0:Z_1:Z_2]\longmapsto
[\overline Z_0:\overline Z_2:\overline Z_1].
$$
Thus $u_0$ is a totally geodesic Lagrangian embedding with boundary on $L$.

In the affine chart $Z_0\ne0$, for $\mathbf w=(w_1,w_2)$ where
$w_j=Z_j/Z_0$, the Fubini--Study metric is
\begin{equation}\label{eq:FS-metric}
g(V,W)=\operatorname{Re}
\frac{(1+|\mathbf w|^2)\langle V,W\rangle_0
 -\langle V,\mathbf w\rangle_0\langle\mathbf w,W\rangle_0}
 {(1+|\mathbf w|^2)^2}.
\end{equation}
The outward conormal of $u_0(\D)$ along its boundary is represented
by $(z,\bar z)$, while the tangent space of $L$ is spanned by
$(iz,0)$ and $(0,i\bar z)$.  Formula~\eqref{eq:FS-metric} shows
that these directions are orthogonal.  Hence $u_0$ is a
free-boundary minimal disk.  Writing $z=re^{i\theta}$ gives
$$
u_0^*g=\frac{2\,dr^2}{(1+2r^2)^2}
       +\frac{2r^2\,d\theta^2}{1+2r^2}.
$$
We equip the disk with the complex structure determined by this
metric and its standard orientation, so that $u_0$ is conformal.
Integration gives $\Area_g(u_0)=2\pi(1-1/\sqrt3)$.

On the other hand, the two disks
$$
h_+(z)=[1:z:1],\quad h_-(z)=[1:1:\bar z]
$$
have boundary on $L$ and are holomorphic and
antiholomorphic, respectively. All three disks lie in the chart $Z_0\ne0$,
identified with $\C^2$, where $L=S^1\times S^1$.
In $\pi_1(L)\cong\Z^2$, the boundary winding classes of
$h_+$, $h_-$, and $u_0$ are $(1,0)$, $(0,-1)$, and $(1,-1)$,
respectively; the minus sign comes from complex conjugation.
Since $\C^2$ is contractible, the boundary map
$\pi_2(\C^2,L)\to\pi_1(L)$ is an isomorphism. Thus
$[u_0]=[h_+]+[h_-]$ in $\pi_2(\C^2,L)$, hence in
$\pi_2(\mathbb{CP}^2,L)$ as well. Each of $h_+$ and $h_-$ has area $\pi/3$
and minimizes area in its relative homotopy class.
A boundary connected sum of $h_+$ and $h_-$, with a connecting
strip of arbitrarily small area, yields
\begin{equation}\label{eq:Clifford-nonmin-area}
\inf_{v\in[u_0]}\Area_g(v)\le\frac{2\pi}{3}
<2\pi\left(1-\frac1{\sqrt3}\right)=\Area_g(u_0).
\end{equation}
Thus the Lagrangian free-boundary minimal disk $u_0$ is not area
minimizing. In fact, it is also unstable.  For small real $t$ put
$$
f_t(z)=\frac{z-t}{1-tz},\quad
 g_t(z)=\frac{z+t}{1+tz},\quad
 F_t(z)=[1:f_t(z):\overline{g_t(z)}].
$$
Because $|f_t|=|g_t|=1$ on $\partial\D$, this is an admissible
variation through maps with boundary on $L$, with $F_0=u_0$.
For $A(t)=\Area_g(F_t)$, differentiation of the area density
gives
$$
A''(0)=8\pi\int_0^1
\frac{r^7+3r^5+11r^3-5r}{(1+2r^2)^{5/2}}\,dr
=\frac{8\pi}{27}(9-8\sqrt3)<0.
$$

Thus the Lagrangian alternative in Theorem~\ref{thm:KE-positive}
can occur, and such disks need not be stable.

\subsection{Sharpness of the Maslov-index criterion}\label{subsec:KE-product-sharpness}
The following family shows that the strict inequality in
Theorem~\ref{thm:intro-A}(ii) is sharp even for immersed area-minimizing disks,
and illustrates the contribution of tangencies at branch points.

Let $M=\mathbb{CP}^1\times\mathbb{CP}^1$ carry the product of
round metrics of Gauss curvature $1$, with complex structure
$J=(J_1,J_2)$ and K\"ahler form $\omega=\omega_1+\omega_2$.
Let $L=S^1\times S^1$ be the product of equators.  Then
$\Ric_M=g$, and $L$ is totally geodesic and Lagrangian.

For integers $p,q\ge1$ and $0\leq a<1$ define
\begin{equation}\label{eq:product-shifted}
T_a(z)=\frac{z-a}{1-az},\quad
u_{p,q,a}(z)=\bigl(z^p,\overline{T_a(z)^q}\bigr).
\end{equation}
Here the unit disk is identified with a hemisphere in each factor.
The map $u_{p,q,a}$ is holomorphic for the parallel complex
structure $J'=(J_1,-J_2)$ and is therefore a conformal minimal
map. Since $|T_a(z)|=1$ for $|z|=1$ and
$$
\frac{zT_a'(z)}{T_a(z)}=\frac{1-a^2}{|z-a|^2}>0
\quad (|z|=1),
$$
the radial conormal is perpendicular to the tangent space of $L$.
Thus these are free-boundary minimal disks.
Writing $u_{p,q,a}=(f,\bar g)$, we have
$$
f'(z)=pz^{p-1},\quad
 g'(z)=\frac{q(1-a^2)(z-a)^{q-1}}{(1-az)^{q+1}}.
$$
For $a>0$, the map $u_{p,q,a}$ is immersed.
The point $0$ is an anticomplex point of multiplicity $p-1$
when $p>1$, and $a$ is a complex point of multiplicity $q-1$
when $q>1$. Hence $R=0$, $Q=p-1$, and $P=q-1$.
The form $\omega' =\omega_1-\omega_2$ calibrates these disks and
vanishes on $L$.  In particular, they minimize area in their
relative homotopy classes, and
\begin{equation}\label{eq:product-shifted-data}
\begin{aligned}
\Area_g(u_{p,q,a})&=2\pi(p+q),\\
\omega(u_{p,q,a})&=2\pi(p-q),\\
\mu_L(u_{p,q,a})&=2(p-q)=2\bigl(Q(u_{p,q,a})-P(u_{p,q,a})\bigr).
\end{aligned}
\end{equation}

Taking $p=m+1$, $q=1$ gives
$P=0$, $Q=m$ and $\mu_L=2m=2Q$ for every $m\ge1$.
These area-minimizing disks are not $J$-holomorphic, so the
strict inequality in Theorem~A(ii) cannot be replaced by a weak
inequality.

When $p=q=1$, both tangency counts, the symplectic area, and the
Maslov index vanish. In the coordinate $z=x+iy$,
\[
(u_{1,1,a}^*\omega)_0=\frac{16a^2}{(1+a^2)^2}\,dx\wedge dy.
\]
The disk $u_{1,1,a}$ is totally real and is Lagrangian if and only if $a=0$.

Finally,
$u_{p,q,0}(z)=(z^p,\bar z^q)$ has a branch point at $0$
when $p,q>1$.  Writing $r=\min\{p-1,q-1\}$, we have
$R=r$, $Q=p-1-r$, and $P=q-1-r$.
For example, $u_{3,2,0}$ has $R=1$, $Q=1$, and $P=0$.
The branch-point convention in the relative complex-point formula
captures the residual tangency at the branch point.

\subsection{A flat minimizing disk of constant nonzero K\"ahler angle}
\label{subsec:flat-constant-angle}
In $\mathbb C^2$ consider the torus
$$
L=\left\{\bigl(\sqrt3(2e^{i\theta}-e^{i\varphi}),
      -2e^{-i\theta}+3e^{-i\varphi}\bigr):
\theta,\varphi\in\mathbb R/(2\pi\mathbb Z)\right\}
$$
and the disk $u(z)=(\sqrt3z,\bar z)$, $|z|\le1$.
Introduce real-linear complex-valued coordinates
$$
\zeta=\frac{\sqrt3w_1+\bar w_2}{2},\quad
\eta=\frac{w_1/\sqrt3+\bar w_2}{2}.
$$
Then $L=\{|\zeta|=2,|\eta|=1\}$ and
$$
\omega_0=2\,dA_\zeta-6\,dA_\eta,
\quad dA_\zeta=d(\operatorname{Re}\zeta)\wedge d(\operatorname{Im}\zeta),
$$
with the analogous notation for $dA_\eta$.  As both area forms
vanish on the tangent planes of $L$, the torus is Lagrangian.
Along the disk, $\zeta\circ u=2z$ and $\eta\circ u=z$, while
$u^*g_0=4|dz|^2$ and $u^*\omega_0=2\,dx\wedge dy$.
Thus $u$ is an immersed conformal minimal disk meeting $L$
orthogonally, with
$$
\Area_{g_0}(u)=4\pi,\quad \omega_0(u)=2\pi,
\quad \cos\alpha=\tfrac12.
$$
In particular, it has no complex or anticomplex points, is not $\pm J_0$-holomorphic, and has $\mu_L(u)=0$.

The closed form $dA_\zeta$ has comass one and vanishes on $L$.
For any disk $v$ relatively homotopic to $u$, the winding number
of $\zeta\circ v|_{\partial\D}$ is one.  Stokes' theorem and
Wirtinger's inequality for the Euclidean projection $\zeta$ give
$$
\Area_{g_0}(v)\ge\int_\D v^*dA_\zeta=4\pi=\Area_{g_0}(u).
$$
Hence area minimization and $\cos\alpha>0$ do not by themselves
force holomorphicity for the given complex structure.
The supporting torus is not minimal, so Theorem~\ref{thm:KE-hyperkahler}
does not apply.

\subsection{Necessity of the Lagrangian boundary hypothesis}
\label{sec:totally-real-example}
Choose the standard orthonormal frame $(e_1,Je_1,e_2,Je_2)$ of
$\mathbb C^2$. Fix $0<c<1$ and $s=\sqrt{1-c^2}$. Let 
$b=cJe_1+s e_2$ and $n=sJe_1-c e_2$, so $e_1,b,n$
are orthonormal. The flat disk $u(x,y)=x e_1+y b$,
$x^2+y^2\le1$, has boundary on the cylinder
$$
L=\left\{\cos\theta\,e_1+\sin\theta\,b+t n:
\theta\in\mathbb R/(2\pi\mathbb Z),\ t\in\mathbb R\right\}.
$$
Along its boundary, the tangent and outward conormal are
$$
T=-\sin\theta\,e_1+\cos\theta\,b,\quad
\eta=\cos\theta\,e_1+\sin\theta\,b,
$$
and $T_{u(e^{i\theta})}L=\operatorname{span}_\mathbb R\{T,n\}$.
The conormal $\eta$ is perpendicular to both $T$ and $n$, so the disk is
free-boundary minimal. It follows from 
$\cos\alpha=\omega(e_1,b)=c\in(0,1)$ that 
$P(u)=Q(u)=0$.

With respect to the standard complex frame,
$$
\det_\mathbb C(T,n)=c\sin\theta-i\cos\theta.
$$
The determinant never vanishes, so $L$ is totally real.
Its phase winds once around the circle, giving $\mu_L(u)=2$.  On the other hand,
$\omega(T,n)=-s\sin\theta\not\equiv 0$,
so $L$ is not Lagrangian.  The disk is not $\pm J$-holomorphic.  Consequently,
$$
\mu_L(u)=2\ne 2\bigl(Q(u)-P(u)\bigr)=0,
\quad \mu_L(u)>2Q(u).
$$
Both the complex-point identity and the holomorphicity criterion
fail if the Lagrangian boundary condition is weakened to total
reality.

\end{document}